# On Construction of Hadamard Matrices

Shipra Kumari[1] and Hrishikesh Mahato[2]

*Department of Mathematics, Central University of Jharkhand, Ranchi-835205, India*
*[1]shipracuj@gmail.com*
*[2]hrishikesh.mahato@cuj.ac.in*

**Abstract.** In this article, a series of Hadamard matrix has been developed using some block matrices with the help of skew Hadamard matrix. Basically an internal structure of skew Hadamard matrix has been changed with some block matrices using kronecker product. For some parameter, Hadamard matrices of order $4t$ where $t$ is an integer, has been found.

## Introduction

Let $H = [h_{ij}]$ be a $(-1, 1)$ square matrix of order $n$, The matrix $H$ is said to be Hadamard matrix of order $n$ if $H^T H = HH^T = nI_n$. It is known that if a Hadamard matrix exists then its order is either 1, or 2 or multiple of 4 [18]. "The existence of Hadamard matrix of order $n = 4t$ with any positive integer $t$" is still an open problem and known as Hadamard conjecture.

Many authors developed several methods to construct Hadamard matrices. There are infinitely many orders of Hadamard matrices have been constructed. But it is also true that there are no constructions have been developed for infinitely many orders (multiple of 4). In 2005 Kharaghani and Tayfeh discovered Hadamard matrix of order 428 (H. Kharaghani, 2005). As of 2008, the unknown smallest order of Hadamard matrix is 668 [20].

The proposed methods of construction of Hadamard are distinguish on the basis of types of an odd integer. The basic idea of the construction is to change the internal structure of skew Hadamard matrix. Basically entries of skew Hadamard matrix are replaced by the block matrices $Q_q = C_q + I_q$, where $C_q$ is a conference matrix, and $J_q$ is a square matrix having all entries 1. The construction requires the existence of skew Hadamard matrix. The Paley construction ensures the existence of skew Hadamard matrices of order $q+1$, where $q \equiv 3 \pmod 4$ is an odd prime. Apart from the Paley constructions there are several methods to construct skew Hadamard matrices, which can be seen in [2, 14].

## Preliminaries

In the main results, the kronecker products[10] of the conference matrix [10] and skew Hadamard matrix[2] has been used in a particular manner. Also it requires the following lemmas and properties in the proof of main results.

We have a construction of conference matrices using extended quadratic character on the field $GF(q)$ explained by following lemma.

**Lemma 2.1.** [10] If the matrix $C_q = [c_{ij}] = \chi(\alpha_j - \alpha_i)$ of order $q$, where $\alpha_i, \alpha_j \in GF(q)$ then

(a) $C_q$ is skew symmetric if $q = 4k+3$ and symmetric if $q = 4k+1$ for $k \in \mathbb{Z}^+$

(b) $C_q C_q^T = qI_q - J_q$

Let $C_q$ be a conference matrix of order $q$, then we define a matrix
$$Q_q = [C_q + I_q] \qquad (1)$$

***Properties 2.1.*** Let $K, L, M, N$ be four matrices then

1. $(K \otimes L)(M \otimes N) = KM \otimes LN$
2. $(K \otimes L) \otimes M = K \otimes (L \otimes M)$

***Properties 2.2.*** Let $Q_q = C_q + 1$ be the matrix defined as equation (1) then

1. $J_q Q_q^T = Q_q^T J_q = J_q$
2. $Q_q J_q = J_q Q_q = J_q$
3. $Q_q(J_q - 2I_q) = (J_q - 2I_q)Q_q$
4. $Q_q^T(J_q - 2I_q) = (J_q - 2I_q)Q_q^T$

**Lemma 2.2..** Let $C_q$ be a circulant matrix and $Q_q = C_q + I_q$, then $QR$ is a symmetric matrix, where $R$ is a back diagonal (0, 1) matrix.

**Lemma 2.3..** Let $R$ be a back diagonal (0, 1) matrix of order $q$ then $R^T R = I_q$.

## Construction

To establish the proposed method skew Hadamard matrix has been used. Consider a skew Hadamard matrix $H_{s+1}$ of order $s+1$ Then $S_{s+1} = H_{s+1} - I_{s+1}$ (2)
is a skew symmetric matrix.

Since
$$H_{s+1} H_{s+1}^T = (s+1)I_{s+1}$$
$$\Rightarrow (S_{s+1} + I_{s+1})(S_{s+1} + I_{s+1})^T = (s+1)I_{s+1}$$
$$\Rightarrow S_{s+1} S_{s+1}^T + S_{s+1} + S_{s+1}^T + I_{s+1} = (s+1)I_{s+1}$$

$$S_{s+1}S_{s+1}^T = sI_{s+1} \qquad [\text{as } S_{s+1}^T = -S_{s+1}] \qquad (3)$$

To develop the methods some block matrices has been used. The concern block matrices defined as follows

Let us define

$$N = K_2 \otimes A + I_2 \otimes B \qquad (4)$$

Where $K_2 = \begin{bmatrix} 0 & 1 \\ -1 & 0 \end{bmatrix}$ and $I_2 = \begin{bmatrix} 1 & 0 \\ 0 & 1 \end{bmatrix}$

Then,

$$NN^T = (K_2 \otimes A + I_2 \otimes B)(K_2 \otimes A + I_2 \otimes B)^T$$

$$= (K_2 \otimes A)(K_2^T \otimes A^T) + (K_2 \otimes A)(I_2 \otimes B^T) + (I_2 \otimes B)(K_2^T \otimes A^T) + (I_2 \otimes B)(I_2 \otimes B^T)$$

$$= K_2 K_2^T \otimes AA^T + K_2 \otimes AB^T + K_2^T \otimes BA^T + I_2 \otimes BB^T$$

$$= I_2 \otimes AA^T + K_2 \otimes AB^T - K_2 \otimes BA^T + I_2 \otimes BB^T$$

$$NN^T = I_2 \otimes AA^T + K_2 \otimes AB^T - K_2 \otimes BA^T + I_2 \otimes BB^T \qquad (5)$$

We have following cases for assumption of $A$ and $B$.

**Case-I :-** If $A = J_q - 2I_q$ and $B = J_q$

From equation (5), we get

$$NN^T = I_2 \otimes (J_q - 2I_q)(J_q - 2I_q)^T + K_2 \otimes (J_q - 2I_q)J_q^T - K_2 \otimes J_q(J_q - 2I_q)^T + I_2 \otimes J_q J_q^T$$

$$= I_2 \otimes (J_q J_q^T - 2J_q - 2J_q^T + 4I_q) + I_2 \otimes J_q J_q^T$$

$$= I_2 \otimes (J_q J_q^T - 4J_q + 4I_q + J_q J_q^T)$$

$$= I_2 \otimes (2qJ_q - 4J_q + 4I_q)$$

So, $\qquad NN^T = I_2 \otimes \left((2q-4)J_q + 4I_q\right) \qquad (6)$

**Case-II :-** If $A = J_q$ and $B = Q_q$, where $Q_q$ is a matrix defined in equation (1).

From equation (5), we get

$$NN^T = I_2 \otimes J_q J_q^T + K_2 \otimes J_q Q_q^T - K_2 \otimes Q_q J_q^T + I_2 \otimes Q_q Q_q^T$$

$$= I_2 \otimes qJ_q + I_2 \otimes \left((q+1)I_q - J_q\right) \quad [\text{using properties 2.3. (1)}]$$

So, $\qquad NN^T = I_2 \otimes \left((q-1)J_q + (q+1)I_q\right) \qquad (7)$

**Case-III:-** If $A = J_q - 2I_q$ and $B = Q_q$, where $Q_q = C_q + I_q$ is a matrix defined in equation (1).

From equation (5), we get

$$NN^T = I_2 \otimes (J_q - 2I_q)(J_q - 2I_q)^T + K_2 \otimes (J_q - 2I_q)Q_q^T - K_2 \otimes Q_q(J_q - 2I_q)^T + I_2 \otimes Q_q Q_q^T$$

$$= I_2 \otimes (J_q J_q^T - 2J_q - 2J_q^T + 4I_q) + I_2 \otimes \left((q+1)I_q - J_q\right)$$

$$= I_2 \otimes (qJ_q - 4J_q + 4I_q) + I_2 \otimes ((q+1)I_q - J_q)$$

$$= I_2 \otimes \left(qJ_q - 4J_q + 4I_q + (q+1)I_q - J_q\right)$$

So, $NN^T = I_2 \otimes ((q-5)J_q + (q+5)I_q)$ (8)

Let us define a $M$ in two different ways as follows

**Case-A :-** For a skew symmetric conference matrix $C_q$ with $q = 4k+1$, for $k \in \mathbb{Z}^+$ define a matrix $M$ as

$$M = \begin{bmatrix} C_q + I_q & C_q + I_q \\ C_q + I_q & -C_q - I_q \end{bmatrix}$$ (9)

Then $MM^T = \begin{bmatrix} C_q + I_q & C_q + I_q \\ C_q + I_q & -C_q - I_q \end{bmatrix} \begin{bmatrix} C_q + I_q & C_q + I_q \\ C_q + I_q & -C_q - I_q \end{bmatrix}^T$

$= \begin{bmatrix} C_q C_q^T + C_q + C_q^T + I_q + C_q C_q^T + C_q + C_q^T + I_q & C_q C_q^T + C_q + C_q^T + I_q - C_q C_q^T - C_q - C_q^T - I_q \\ C_q C_q^T + C_q + C_q^T + I_q - C_q C_q^T - C_q - C_q^T - I_q & C_q C_q^T + C_q + C_q^T + I_q + C_q C_q^T + C_q + C_q^T + I_q \end{bmatrix}$

$= \begin{bmatrix} 2(C_q C_q^T + I_q) & 0 \\ 0 & 2(C_q C_q^T + I_q) \end{bmatrix}$

$= \begin{bmatrix} 2(q+1)I_q - 2J_q & 0 \\ 0 & 2(q+1)I_q - 2J_q \end{bmatrix} = I_2 \otimes (2(q+1)I_q - 2J_q)$ (10)

**Case-B** For a symmetric conference matrix $C_q$ with $q \equiv 1 \pmod 4$ define a matrix $M$ as

$$M = \begin{bmatrix} C_q + I_q & -C_q + I_q \\ -C_q + I_q & -C_q - I_q \end{bmatrix}$$ (11)

Now,

$MM^T = \begin{bmatrix} C_q + I_q & -C_q + I_q \\ -C_q + I_q & -C_q - I_q \end{bmatrix} \begin{bmatrix} C_q + I_q & -C_q + I_q \\ -C_q + I_q & -C_q - I_q \end{bmatrix}^T$

$= \begin{bmatrix} (C_q + I_q)(C_q^T + I_q) + (-C_q + I_q)(-C_q^T + I_q) & (C_q + I_q)(-C_q^T + I_q) + (-C_q + I_q)(-C_q^T - I_q) \\ (-C_q + I_q)(C_q^T + I_q) + (-C_q - I_q)(-C_q^T + I_q) & (-C_q + I_q)(-C_q^T + I_q) + (-C_q - I_q)(-C_q^T - I_q) \end{bmatrix}$

$= \begin{bmatrix} C_q C_q^T + C_q + C_q^T + I_q + C_q C_q^T - C_q - C_q^T + I_q & -C_q C_q^T + C_q - C_q^T + I_q + C_q C_q^T + C_q - C_q^T - I_q \\ -C_q C_q^T - C_q + C_q^T + I_q + C_q C_q^T - C_q + C_q^T - I_q & C_q C_q^T - C_q - C_q^T + I_q + C_q C_q^T + C_q + C_q^T + I_q \end{bmatrix}$

$= \begin{bmatrix} 2(C_q C_q^T + I_q) & 0 \\ 0 & 2(C_q C_q^T + I_q) \end{bmatrix}$

$= \begin{bmatrix} 2(q+1)I_q - 2J_q & 0 \\ 0 & 2(q+1)I_q - 2J_q \end{bmatrix} = I_2 \otimes (2(q+1)I_q - 2J_q)$ (12)

**Theorem 3.1.** Let $q$ be any odd integer such that $q \equiv 1 \pmod{4}$ and $s$ is an odd integer of the form of $3 \pmod{4}$ with $s = q - 2$. If there exist a skew Hadamard matrix of order $s+1$ and a circulant conference matrix of order $q$ then there is a Hadamard matrix of order $2q(s+1)$.

*Proof:* Consider a circulant conference matrix $C_q$ then either $C_q$ is symmetric or $C_q R$ is symmetric. Without loss of generality, consider $C_q$ is symmetric.

Let $H_{s+1}$ be a skew Hadamard matrix, so
$$H_{s+1} = S_{s+1} + I_{s+1}$$

Let us define
$$H' = S_{s+1} \otimes M + I_{s+1} \otimes N$$

Where $M = \begin{bmatrix} C_q + I_q & -C_q + I_q \\ -C_q + I_q & -C_q - I_q \end{bmatrix}$ defined in equation (11) and $N = \begin{bmatrix} J_q - 2I_q & J_q \\ -J_q & J_q - 2I_q \end{bmatrix}$ is defined in equation case-I.

Using properties 2.3, It can be easily verified that $MN^T = NM^T$.
Then we see that $H'$ forms a Hadamard matrix.

We have
$$H'H'^T = (S_{s+1} \otimes M + I_{s+1} \otimes N)(S_{s+1} \otimes M + I_{s+1} \otimes N)^T$$
$$= (S_{s+1} \otimes M)(S_{s+1}^T \otimes M^T) + (S_{s+1} \otimes M)(I_{s+1} \otimes N^T) + (I_{s+1} \otimes N)(S_{s+1}^T \otimes M^T) + (I_{s+1} \otimes N)(I_{s+1} \otimes N^T)$$
$$= S_{s+1}S_{s+1}^T \otimes MM^T + S_{s+1} \otimes MN^T + S_{s+1}^T \otimes NM^T + I_{s+1} \otimes NN^T$$

$$H'H'^T = S_{s+1}S_{s+1}^T \otimes MM^T + I_{s+1} \otimes NN^T$$

Now substituting the value of $S_{s+1}S_{s+1}^T$ from equation (3), that of $NN^T$ from equation (6) and that of $MM^T$ from equation (12), we have
$$H'H'^T = sI_{s+1} \otimes (I_2 \otimes (2(q+1)I_q - 2J_q) + I_{s+1} \otimes (I_2 \otimes ((2q-4)J_q + 4I_q))$$
$$= sI_{s+1} \otimes (I_2 \otimes 2(q+1)I_q - I_2 \otimes 2J_q) + I_{s+1} \otimes (I_2 \otimes (2q-4)J_q + I_2 \otimes 4I_q)$$
$$= sI_{s+1} \otimes (2(q+1)I_{2q} - 2I_2 \otimes J_q) + I_{s+1} \otimes ((2q-4)I_2 \otimes J_q + 4I_{2q})$$
$$= sI_{s+1} \otimes 2(q+1)I_{2q} - sI_{s+1} \otimes 2I_2 \otimes J_q + I_{s+1} \otimes (2q-4)I_2 \otimes J_q + I_{s+1} \otimes 4I_{2q}$$
$$= 2s(q+1)I_{2q(s+1)} - 2sI_{2(s+1)} \otimes J_q + (2q-4)I_{2(s+1)} \otimes J_q + 4I_{2q(s+1)}$$
$$= (2s(q+1) + 4)I_{2q(s+1)} + (2q - 4 - 2s)I_{2(s+1)} \otimes J_q$$

Since, it is given that $s = q - 2 \Rightarrow 2s = 2q - 4$ and also $2s + 4 = 2q$.

Thus
$$H'H'^T = (2q + 2sq)I_{2q(s+1)}$$
i.e. $\quad H'H'^T = 2q(s+1)I_{2q(s+1)}$

Thus we obtain a Hadamard matrix $H'$ of order $2q(s+1)$. $\square$

**Corollary 3.1.** There exist a Hadamard matrix of order $2q(s+1)$, where $q = u^v \equiv 1 \pmod{4}$ and $s = q - 2 = p^r$, and both $p$ & $u$ are odd primes and $r, v \in \mathbb{Z}^+$.

*Proof.* Since $s = p^r$, using Paley construction-I there exist a skew Hadamard matrix of order $s+1$. Using lemma 2.1, a conference matrix of order $q$ may be constructed. Thus there exist a Hadamard matrix of order $2q(s+1)$.

**Note:-** If $H_n = S_n + I_n$ is a skew Hadamard matrix, then

$$H_{2n} = \begin{bmatrix} S_n + I_n & S_n + I_n \\ S_n - I_n & -S_n + I_n \end{bmatrix}$$

is also a skew Hadamard matrix [2].

**Theorem 3.2.** There exist a Hadamard matrix of order $2q(s+1)$ if a skew Hadamard matrix of order $s+1$ and a conference matrix of order $q$ exist, where $s \equiv \dfrac{q-1}{2} \equiv 3 \pmod{4}$.

*Proof.:-* Let $H_{s+1}$ be a skew Hadamard matrix, so

$$H_{s+1} = S_{s+1} + I_{s+1}$$

Let us define

$$H' = S_{s+1} \otimes M + I_{s+1} \otimes N$$

Where $M = \begin{bmatrix} C_q + I_q & C_q + I_q \\ C_q + I_q & -C_q - I_q \end{bmatrix}$ is defined in equation (9) and $N = \begin{bmatrix} J_q & Q_q \\ -Q_q & J_q \end{bmatrix}$ is defined in Case-II. Using properties 2.3, it can be easily verified that $MN^T = NM^T$.

We see that $H'$ forms a Hadamard matrix.
We have

$$H' H'^T = (S_{s+1} \otimes M + I_{s+1} \otimes N)(S_{s+1} \otimes M + I_{s+1} \otimes N)^T$$

$$= (S_{s+1} \otimes A')(S^T_{s+1} \otimes A'^T) + (S_{s+1} \otimes M)(I_{s+1} \otimes N^T) + (I_{s+1} \otimes N)(S^T_{s+1} \otimes M^T) + (I_{s+1} \otimes N)(I_{s+1} \otimes N^T)$$

$$= S_{s+1} S^T_{s+1} \otimes MM^T + S_{s+1} \otimes MN^T + S^T_{s+1} \otimes NM^T + I_{s+1} \otimes NN^T$$

$$= S_{s+1} S^T_{s+1} \otimes MM^T + I_{s+1} \otimes NN^T$$

$$= sI_{s+1} \otimes \left( I_2 \otimes (2(q+1)I_q - 2J_q) \right) + I_{s+1} \otimes \left( I_2 \otimes ((q+1)I_q + (q-1)J_q) \right)$$

$$= sI_{s+1} \otimes \left( 2(q+1)I_{2q} - 2I_2 \otimes J_q \right) + I_{s+1} \otimes \left( (q+1)I_{2q} + (q-1)I_2 \otimes J_q \right)$$

$$= sI_{s+1} \otimes 2(q+1)I_{2q} - sI_{s+1} \otimes 2I_2 \otimes J_q + I_{s+1} \otimes (q+1)I_{2q} + I_{s+1} \otimes (q-1)I_2 \otimes J_q$$

$$= (2sq + 2s + q + 1)I_{2q(s+1)} + (q - 1 - 2s)I_{2(s+1)} \otimes J_q$$

Since $2s = q - 1 \Rightarrow q - 1 - 2s = 0$, and $2s + 1 = q$ So

$$H'H'^T = 2q(s+1)I_{2q(s+1)}.$$

This shows that $H'$ is a Hadamard matrix of order $2q(s+1)$.

**Corollary 3.2.** There exist a Hadamard matrix of order $2q(s+1)$, where $s = p^r \equiv 3 \pmod 4$ and $q = 2s+1 = u^v$, both $p$ & $u$ are odd primes and $r, v \in \mathbb{Z}^+$.

*Proof.* Here $s = p^r$, using Paley construction-I there exist a skew Hadamard matrix of order $s+1$. Using lemma 2.1, a conference matrix of order $q$ may be constructed. Thus there exist a Hadamard matrix of order $2q(s+1)$.

**Remark 3.1.** For twin primes $p$ and $r$ we have a conference matrix of order $q = pr$ [10].

**Theorem 3.3.** There exist a Hadamard matrix of order $2q(s+1)$ if a skew Hadamard matrix of order $s+1$ and a conference matrix of order $q$ exist, where $s \equiv 3 \pmod 4$ and $q \equiv 3 \pmod 8$ with $2s = q-5$.

*Proof.* Consider a skew Hadamard matrix $H_{s+1}$. Then $H_{s+1}$
$$H_{s+1} = S_{s+1} + I_{s+1}$$

Let us define
$$H' = S_{s+1} \otimes M + I_{s+1} \otimes N$$

Where $M = \begin{bmatrix} C_q + I_q & C_q + I_q \\ C_q + I_q & -C_q - I_q \end{bmatrix}$ is defined in equation (9) and $N = \begin{bmatrix} J_q - 2I_q & Q_q \\ -Q_q & J_q - 2I_q \end{bmatrix}$ is defined in Case-III. Using properties 2.3, it can be easily verified that $MN^T = NM^T$.

We see that $H'$ form a Hadamard matrix of order $2q(s+1)$.
We have,
$$H'H'^T = (S_{s+1} \otimes M + I_{s+1} \otimes N)(S_{s+1} \otimes M + I_{s+1} \otimes N)^T$$
$$= (S_{s+1} \otimes M)(S_{s+1}^T \otimes M^T) + (S_{s+1} \otimes M)(I_{s+1} \otimes N^T) + (I_{s+1} \otimes N)(S_{s+1}^T \otimes M^T) + (I_{s+1} \otimes N)(I_{s+1} \otimes N^T)$$
$$= S_{s+1}S_{s+1}^T \otimes MM^T + S_{s+1} \otimes MN^T + S_{s+1}^T \otimes NM^T + I_{s+1} \otimes NN^T$$
$$= S_{s+1}S_{s+1}^T \otimes MM^T + I_{s+1} \otimes NN^T$$
$$= sI_{s+1} \otimes \left(I_2 \otimes (2(q+1)I_q - 2J_q)\right) + I_{s+1} \otimes \left(I_2 \otimes ((q-5)J_q + (q+5)I_q)\right)$$
$$= sI_{s+1} \otimes \left(I_2 \otimes 2(q+1)I_q - I_2 \otimes 2J_q\right) + I_{s+1} \otimes \left(I_2 \otimes (q-5)J_q + I_2 \otimes (q+5)I_q\right)$$
$$= sI_{s+1} \otimes \left(2(q+1)I_{2q} - 2I_2 \otimes J_q\right) + I_{s+1} \otimes \left((q-5)I_2 \otimes J_q + (q+5)I_{2q}\right)$$
$$= sI_{s+1} \otimes 2(q+1)I_{2q} - sI_{s+1} \otimes 2I_2 \otimes J_q + I_{s+1} \otimes (q-5)I_2 \otimes J_q + I_{s+1} \otimes (q+5)I_{2q}$$
$$= (2s(q+1) + (q+5))I_{2q(s+1)} + (q-5-2s)I_{2(s+1)} \otimes J_q$$
Since $2s = q-5 \Rightarrow 2s+5 = q$, and $q-5-2s = 0$ So

$$H'H'^T = 2q(s+1)I_{2q(s+1)}.$$

This shows that $H'$ is a Hadamard matrix of order $2q(s+1)$.

**Theorem 3.4.** There exist a Hadamard matrix of order $2q(s+1)$ if a skew symmetric Hadamard matrix of order $s+1$ and conference matrix of order $q$ exist, where $s = q - 4 \equiv 3 \pmod 4$.

*Proof.* Consider a skew Hadamard matrix $H_{s+1}$. Then $H_{s+1} = S_{s+1} + I_{s+1}$.

Let us define
$$H' = S_{s+1} \otimes M + I_{s+1} \otimes N$$

Where $M = \begin{bmatrix} C_q + I_q & C_q + I_q \\ C_q + I_q & -C_q - I_q \end{bmatrix}$ is defined in equation (9) and $N$ is defined as

$$N = \begin{bmatrix} J_q - 2I_q & J_q - 2I_q \\ J_q - 2I_q & -(J_q - 2I_q) \end{bmatrix}$$

Using properties 2.3, it can be easily verified that $MN^T = NM^T$.

We see that $H'$ forms a Hadamard matrix.

Firstly we will calculate $NN^T$. We have

$$NN^T = \begin{bmatrix} J_q - 2I_q & J_q - 2I_q \\ J_q - 2I_q & -(J_q - 2I_q) \end{bmatrix} \begin{bmatrix} J_q - 2I_q & J_q - 2I_q \\ J_q - 2I_q & -(J_q - 2I_q) \end{bmatrix}^T$$

$$= \begin{bmatrix} J_q - 2I_q & J_q - 2I_q \\ J_q - 2I_q & -(J_q - 2I_q) \end{bmatrix} \begin{bmatrix} J_q^T - 2I_q & J_q^T - 2I_q \\ J_q^T - 2I_q & -(J_q^T - 2I_q) \end{bmatrix}$$

$$= \begin{bmatrix} (J_q - 2I_q)(J_q^T - 2I_q) + (J_q - 2I_q)(J_q^T - 2I_q) & (J_q - 2I_q)(J_q^T - 2I_q) - (J_q - 2I_q)(J_q^T - 2I_q) \\ (J_q - 2I_q)(J_q^T - 2I_q) - (J_q - 2I_q)(J_q^T - 2I_q) & (J_q - 2I_q)(J_q^T - 2I_q) + (J_q - 2I_q)(J_q^T - 2I_q) \end{bmatrix}$$

$$= \begin{bmatrix} J_q J_q^T - 2J_q - 2J_q^T + 4I_q + J_q J_q^T - 2J_q - 2J_q^T + 4I_q & J_q J_q^T - 2J_q - 2J_q^T + 4I_q - J_q J_q^T + 2J_q + 2J_q^T - 4I_q \\ J_q J_q^T - 2J_q - 2J_q^T + 4I_q - J_q J_q^T + 2J_q + 2J_q^T - 4I_q & J_q J_q^T - 2J_q - 2J_q^T + 4I_q + J_q J_q^T - 2J_q - 2J_q^T + 4I_q \end{bmatrix}$$

$$= \begin{bmatrix} 2(J_q J_q^T - 2J_q - 2J_q^T + 4I_q) & 0 \\ 0 & 2(J_q J_q^T - 2J_q - 2J_q^T + 4I_q) \end{bmatrix}$$

$$= \begin{bmatrix} 2(qJ_q - 4J_q + 4I_q) & 0 \\ 0 & 2(qJ_q - 4J_q + 4I_q) \end{bmatrix} \qquad \text{As } J_q^T = J_q \text{ and } J_q J_q^T = qJ_q$$

$$= \begin{bmatrix} 2(q-4)J_q + 8I_q & 0 \\ 0 & 2(q-4)J_q + 8I_q \end{bmatrix}$$

So, $$NN^T = I_2 \otimes \left(2(q-4)J_q + 8I_q\right) \qquad (13)$$

Now, we have
$$H'H'^T = (S_{s+1} \otimes M + I_{s+1} \otimes N)(S_{s+1} \otimes M + I_{s+1} \otimes N)^T$$

$$\begin{aligned}
&= \left(S_{s+1} \otimes M + I_{s+1} \otimes N\right)\left(S_{s+1}^T \otimes M^T + I_{s+1} \otimes N^T\right) \\
&= (S_{s+1} \otimes M)(S_{s+1}^T \otimes M^T) + (S_{s+1} \otimes M)(I_{s+1} \otimes N^T) + (I_{s+1} \otimes N)(S_{s+1}^T \otimes M^T) + (I_{s+1} \otimes N)(I_{s+1} \otimes N^T) \\
&= S_{s+1}S_{s+1}^T \otimes MM^T + I_{s+1} \otimes NN^T \\
&= sI_{s+1} \otimes \left(I_2 \otimes (2(q+1)I_q - 2J_q)\right) + I_{s+1} \otimes \left(I_2 \otimes ((2q-8)J_q + 8I_q)\right) \\
&= sI_{s+1} \otimes \left(I_2 \otimes 2(q+1)I_q - I_2 \otimes 2J_q\right) + I_{s+1} \otimes \left(I_2 \otimes (2q-8)J_q + I_2 \otimes 8I_q\right) \\
&= sI_{s+1} \otimes \left(2(q+1)I_{2q} - 2I_2 \otimes J_q\right) + I_{s+1} \otimes \left((2q-8)I_2 \otimes J_q + 8I_{2q}\right) \\
&= 2s(q+1)I_{2q(s+1)} - 2sI_{2(s+1)} \otimes J_q + (2q-8)I_{2(s+1)} \otimes J_q + 8I_{2q(s+1)} \\
&= \left(2s(q+1) + 8\right)I_{2q(s+1)} + (2q-8-2s)I_{2(s+1)} \otimes J_q
\end{aligned}$$

It has been considered that $s = q - 4 \Rightarrow 2s = 2q - 8$ and $2s + 8 = 2q$

So, $H'H'^T = 2q(s+1)I_{2q(s+1)}$

This shows that $H'$ is a Hadamard matrix of order $2q(s+1)$.

## CONCLUSION

In proposed methods the known skew Hadamard matrices of smaller order have been used to construct a series of Hadamard matrices. Basically we have replaced the entries of skew Hadamard matrix with the block matrices $Q_q = C_q + I_q$ and $J_q$. In this manner, we get a series of Hadamard matrices.

## REFERENCES

1. A. Hedayat and W. D. Wallis, "Hadamard Matrices and Their Applications" The Annals of Statistics, $6$, pp $1184 - 1238 \, (1978)$.
2. C. Koukouvinos, S. Stylianou, "On skew-Hadamard matrices", Discrete Mathematics, $308$, pp $2723 - 2731 (2008)$.
3. H. Kharaghani and B. Tayfeh-Rezaie, "A Hadamard matrix of order 428", Journal of Combinatorial Designs, $13$, pp $435 - 440 (2005)$.
4. H. Mahato, "On existence of williamson symmetric circulant matrices", Bulletin of the Institute of Mathematics Academia Sinica (New Series) Vol $6$, No. 1, pp $27 - 39 (2011)$.
5. J. M. Goethals, and J. J. Seidel, "Orthogonal Matrices with Zero Diagonal", Geometrics and Combinatorics, pp $257 - 266 (1991)$.
6. J. Seberry and M. Yamada, "Hadamard matrices, sequences and block designs in Contemporary Design Theory: A Collection of Surveys", John Wiley and Sons, pp $431 - 560 (1992)$
7. J. Steepleton, "Construction of Hadamard Matrices", TRACE: Tennessee Research and Creative Exchange, $5 (2019)$.
8. K. J. Horadam, Hadamard Matrices and Their Application, Princeton University Press, 2007.

9. L. D. Baumert and M. Hall, Jr., "A new construction for Hadamard matrices", Bulletin of American Mathematical Society, $71$, pp $169-170(1965)$.
10. Marshall, Jr. Hall, "Combinatorial Theory", Wiley Interscience Series in Discrete Mathematics, MAc-Wiley, $1986,$ second edition.
11. M. Mitrouli, "Sylvester Hadamard matrices revisited", Special Matrices, $2$, pp $120-124(2014)$.
12. M. Mitrouli, "A sign test for detecting the equivalence of Sylvester Hadamard matrices" Numerical Algorithm, $57$, pp $169-186(2011)$.
13. M. Miyamoto, "A Construction of Hadamard Matrices", Journal of Combinatorial Theory, Series A $57$, pp $86-108(1991)$.
14. N. A. Malik and H. Mahato, "A New Approach to Good Matrices", The Mathematical Student, 90, pp 161-175 (2020).
15. P. Delsarte, J. M. Goethals and J. J. Seidel, "Orthogonal Matrices with Zero Diagonal-II", Canadian Journal of Mathematics, $5$ $5$ pp $816-832(1971)$.
16. P.K.Manjhi and A. Kumar, "ON THE CONSTRUCTION OF HADAMARD MATRICES", International Journal of Pure and Applied Mathematics, Vol $120$, No. 1, pp $51-58(2018)$.
17. R. E. A. C. Paley, "On orthogonal matrices", Journal of Mathematical Physics, $12$, pp $311-320(1933)$.
18. S. Kumari, H. Mahato, " Extension of Paley Construction for Hadamard Matrices", Journal of Interdisciplinary Cycle Research, Vol $XIII$, Issue III pp $1506-1525(2021)$
19. https://arxiv.org/abs/1912.10757v1
20. W. V. Nishadi, A. A. Perera, "On Construction of Hadamard Matrices", Iconic Research and Engineering Journals", vol $3$, No. $11,$ pp $261-264(2020)$.